\documentclass[onecolumn,11pt]{article}
\usepackage[top=1in, bottom=1in, left=1in, right=1in]{geometry}
\usepackage{amsmath, amsthm, amssymb}
\usepackage{graphicx}
\usepackage{epstopdf}
\usepackage{overpic}
\usepackage{cancel}
\usepackage{rotating}
\usepackage{url}
\usepackage{caption}
\usepackage{color}
\usepackage{rotating}
\usepackage{multirow}
\usepackage{wrapfig}
\usepackage{mathtools}
\usepackage{arydshln}
\usepackage{setspace}
\usepackage{palatino} 
\usepackage[bottom,flushmargin,hang,multiple]{footmisc}
\usepackage{lipsum}
\newcommand\blfootnote[1]{%
  \begingroup
  \renewcommand\thefootnote{}\footnote{#1}%
  \addtocounter{footnote}{-1}%
  \endgroup
}

\newcommand{\bx}{\mathbf{x}}

\newcommand{\by}{\mathbf{y}}

\newcommand{\bX}{\mathbf{X}}

\newcommand{\bXi}{\boldsymbol{\Xi}}
\newcommand{\bxi}{\boldsymbol{\xi}}

\newcommand{\bTheta}{\boldsymbol{\Theta}}

\usepackage{textcomp,xcolor}
\definecolor{MatlabCellColour}{RGB}{250,250,250}
\definecolor{MatPurp}{rgb}{.625,.1406,.9375}
\usepackage{listings}

\lstdefinestyle{customc}{
  belowcaptionskip=.25\baselineskip,
  breaklines=true,
  frame=L,
  xleftmargin=\parindent,
  language=Matlab,
  showstringspaces=false,
  basicstyle=\small\ttfamily,
  keywordstyle=\bfseries\color{white!30!black},
  identifierstyle=\color{blue},  
  commentstyle=\itshape\color{green!60!black},
  stringstyle=\color{MatPurp},
  backgroundcolor=\color{MatlabCellColour}
 }
 
 \usepackage{tikz}
\usetikzlibrary{positioning,shapes,shadows,arrows,calc}

\tikzstyle{block} = [draw,line width=1.2pt, fill=white!30, rectangle, 
    minimum height=3em, minimum width=3em, rounded corners]    
\tikzstyle{blocksys} = [draw, line width=1.2pt,fill=blue!20, rectangle, 
    minimum height=5.em, minimum width=3em, rounded corners]        
\tikzstyle{blockest} = [draw, line width=1.2pt,fill=yellow!40, rectangle, 
    minimum height=5.em, minimum width=3em, rounded corners]        
\tikzstyle{blocky} = [draw, line width=1.2pt,fill=yellow!40, rectangle, 
    minimum height=2.5em, minimum width=2.5em, rounded corners]    
\tikzstyle{blockp} = [draw,line width=1.2pt, fill=purple!30, rectangle, 
    minimum height=3em, minimum width=3em, rounded corners]       
\tikzstyle{blockb} = [draw,line width=1.2pt, fill=blue!20, rectangle, 
    minimum height=2.5em, minimum width=2.5em, rounded corners]     
\tikzstyle{blockr} = [draw, line width=1.2pt,fill=red!30, rectangle, 
    minimum height=2.5em, minimum width=2.5em, rounded corners]                 
    \tikzstyle{blockgr} = [draw, line width=1.2pt,fill=black!10, rectangle, 
    minimum height=2.5em, minimum width=2.5em, rounded corners]                  

\tikzstyle{blockg} = [draw, line width=1.2pt,fill=green!20, rectangle, 
    minimum height=3em, minimum width=3em, rounded corners]    
\tikzstyle{sum} = [draw,line width=1.2pt, fill=white!10, circle, node distance=1cm]
\tikzstyle{input} = [coordinate]
\tikzstyle{output} = [coordinate]
\tikzstyle{pinstyle} = [pin edge={to-,thin,black}]
\tikzstyle{pinstyle2} = [pin edge={to-,thin,black}]
\tikzstyle{line}=[-, line width=1.2pt]
\tikzstyle{blockestbig} = [draw, line width=1.2pt,fill=yellow!40, rectangle,text centered,  
    minimum height=6.em,  text width=9.em, rounded corners]    
\tikzstyle{blocky2} = [draw, line width=1.2pt,fill=yellow!40, rectangle,text centered,  
    minimum height=5.em,  text width=4.em, rounded corners]    
    \tikzstyle{blocksysbig} = [draw, line width=1.2pt,fill=blue!20, rectangle, 
        text centered,minimum height=6.em, text width=9.25em, rounded corners]    
\tikzstyle{blockK} = [draw, line width=1.2pt,fill=yellow!40, rectangle,text centered,  
    minimum height=3.5em,  text width=5em, rounded corners]               
    \tikzstyle{blocky3} = [draw, line width=1.2pt,fill=yellow!40, rectangle, 
        text centered,minimum height=6.em, text width=9.25em, rounded corners]        
  \tikzstyle{blockr2} = [draw, line width=1.2pt,fill=red!30, rectangle,text centered,  
    minimum height=3.em,  text width=4.5em, rounded corners]          

    \tikzstyle{blocky4} = [draw, line width=1.2pt,fill=yellow!40, rectangle, 
        text centered,minimum height=4.5em, text width=6.5em, rounded corners]        
    \tikzstyle{blockb4} = [draw, line width=1.2pt,fill=blue!20, rectangle, 
        text centered,minimum height=4.5em, text width=6.5em, rounded corners]

\definecolor{blue}{rgb}{0,0,1}
\definecolor{darkgreen}{rgb}{0,0.5,0}
\definecolor{red}{rgb}{1,0,0}
\definecolor{teal}{rgb}{0,0.5,0.7}

\DeclareMathOperator*{\argmin}{arg\rm{}min}

\newenvironment{lbmatrix}[1]
  {\left[\array{@{}*{#1}{c}@{}}}
  {\endarray\right]}

\title{\LARGE{\vspace{-.35in}Koopman invariant subspaces and finite linear representations of nonlinear dynamical systems for control\vspace{-.05in}}}
\author{Steven L. Brunton$^{1*}$, Bingni W. Brunton$^{2}$, Joshua L. Proctor$^3$, J. Nathan Kutz$^4$\\
\footnotesize{$^1$ Department of Mechanical Engineering, University of Washington, Seattle, WA 98195, United States}\\
\footnotesize{$^2$ Department of Biology, University of Washington, Seattle, WA 98195, United States}\\
\footnotesize{$^3$Institute for Disease Modeling, Bellevue, WA 98004, United States}\\ 
\footnotesize{$^4$ Department of Applied Mathematics, University of Washington, Seattle, WA 98195, United States}\vspace{-.35in}
}
\date{}
\begin{document}
\maketitle
\blfootnote{$^*$ Corresponding author Email: sbrunton@uw.edu (S.L. Brunton).}
\begin{abstract}
In this work, we explore finite-dimensional linear representations of nonlinear dynamical systems by restricting the Koopman operator to an invariant subspace spanned by specially chosen observable functions.  
The Koopman operator is an infinite-dimensional linear operator that evolves functions of the state of a dynamical system \cite[Koopman 1931]{Koopman1931pnas}.  
Dominant terms in the Koopman expansion are typically computed using dynamic mode decomposition (DMD).  
DMD uses linear measurements of the state variables, and it has recently been shown that this may be too restrictive for nonlinear systems~\cite[Williams et al., 2015]{Williams2015jnls}.  
Choosing the right \emph{nonlinear} observable functions to form an invariant subspace where it is possible to obtain linear reduced-order models, especially those that are useful for control, is an open challenge.

Here, we investigate the choice of observable functions for Koopman analysis that enable the use of optimal linear control techniques on nonlinear problems.  
First, to include a cost on the state of the system, as in linear quadratic regulator (LQR) control, it is helpful to include these states in the observable subspace, as in DMD. 
However, we find that this is only possible when there is a single isolated fixed point, as systems with multiple fixed points or more complicated attractors are not globally topologically conjugate to a finite-dimensional linear system, and cannot be represented by a finite-dimensional linear Koopman subspace that includes the state. 
We then present a data-driven strategy to identify relevant observable functions for Koopman analysis by leveraging a new algorithm to determine relevant terms in a dynamical system by $\ell_1$-regularized regression of the data in a nonlinear function space~\cite[Brunton et al., 2015]{Brunton2015arxiv}; we also show how this algorithm is related to DMD.
Finally, we demonstrate the usefulness of nonlinear observable subspaces in the design of Koopman operator optimal control laws for fully nonlinear systems using techniques from linear  optimal control.  
\\

\end{abstract}

\section{Introduction}\label{sec:introduction}
Koopman spectral analysis provides an operator-theoretic perspective to dynamical systems, which complements the more standard geometric~\cite{guckenheimer_holmes} and probabilistic perspectives.  
In the early 1930s~\cite{Koopman1931pnas,Koopman1932pnas}, B. O. Koopman showed that nonlinear dynamical systems associated with Hamiltonian flows could be analyzed with an infinite dimensional linear operator on the Hilbert space of observable functions.  
For Hamiltonian fluids, the Koopman operator is unitary, meaning that the inner product of any two observable functions remains unchanged by the operator.  
Unitarity is a familiar concept, as the discrete Fourier transform (DFT) and the proper orthogonal decomposition (POD)~\cite{HLBR_turb} both provide unitary coordinate transformations.  
In the original paper~\cite{Koopman1931pnas}, Koopman drew connections between the Koopman eigenvalue spectrum and conserved quantities, integrability, and ergodicity.  
Recently, it was shown that level sets of the Koopman eigenfunctions form invariant partitions of the state-space of a dynamical system~\cite{Budivsic2012physd}; in particular, eigenfunctions of the Koopman operator may be used to analyze the ergodic partition~\cite{mezic1999chaos,Budivsic2009cdc}.  
Koopman analysis has also been recently shown to generalize the Hartman-Grobman theorem to the entire basin of attraction of a stable or unstable equilibrium point or periodic orbit~\cite{Lan2013physd}.  
For more information there are a number of excellent in-depth reviews on Koopman analysis by Mezi\'c et al.~\cite{Budivsic2012chaos,Mezic2013arfm}.  

Koopman analysis has been at the focus of recent data-driven efforts to characterize complex systems, since the work of Mezi\'c and Banaszuk~\cite{Mezic2004physicad} and Mezi\'c~\cite{Mezic2005nd}. 
There is considerable interest in obtaining finite-rank approximations to the linear Koopman operator that propagate the original nonlinear dynamics.    
This is especially promising for the potential control of nonlinear systems~\cite{Proctor2014arxiv}.  
However, by introducing the Koopman operator, we trade nonlinear dynamics for infinite-dimensional linear dynamics, introducing new challenges. 
Finite-dimensional linear approximations of the Koopman operator may be useful to model the dynamics on an attractor, and those that explicitly advance the state may also be useful for control. 
Any set of Koopman eigenfunctions will form a Koopman-invariant subspace, resulting in an exact finite-dimensional linear model.  
Unfortunately, many dynamical systems do not admit a finite-dimensional Koopman-invariant subspace that also spans the state; in fact, this is only possible for systems with an isolated fixed point.
It may be possible to recover the state from the Koopman eigenfunctions, but determining the eigenfunctions and inverting for the state may both be challenging.  

Dynamic mode decomposition (DMD), introduced in the fluid dynamics community~\cite{Schmid2008aps,Rowley2009jfm,Schmid2010jfm,Tu2014jcd}, provides a practical numerical framework for Koopman mode decomposition.  
DMD implicitly uses linear observable functions, such as direct velocity field measurements from particle image velocimetry (PIV).  
In other words, the observable function is an identity map on the fluid flow state.  
This set of linear observables is too limited to describe the rich dynamics observed in fluids or other nonlinear systems.  
Recently, DMD has been extended to include a richer set of nonlinear observable functions, providing the ability to effectively analyze nonlinear systems~\cite{Williams2015jnls}.  
Because of the extreme cost associated with this extended DMD for high-dimensional systems, a variation using the kernel trick from machine learning has been implemented to make the cost of extended DMD equivalent to traditional DMD, but retaining the benefit of nonlinear observables~\cite{Williams2014arxivA}.  
However, choosing the correct nonlinear observable functions to use for a given system, and how they will impact the performance of Koopman mode decomposition and reduction, is still an open problem.  
Presently, these observable functions are either determined using information about the right-hand side of the dynamics (i.e., knowing that the Navier-Stokes equations have quadratic nonlinearities, etc.) or by brute-force trial and error in a particular basis for Hilbert space (i.e., trying many different polynomial functions).

In this work, we explore the identification of observable functions that span a finite-dimensional subspace of Hilbert space which remains invariant under the Koopman operator (i.e., a Koopman-invariant subspace spanned by eigenfunctions of the Koopman operator). 
When this subspace includes the original states, we obtain a finite-dimensional linear dynamical system on this subspace that also advances the original state directly.  
We utilize a new algorithm, the sparse identification of nonlinear dynamics (SINDy)~\cite{Brunton2015arxiv}, to first identify the right-hand side dynamics of the nonlinear system.  
Next, we choose observable functions such that these dynamics are in the span.  
Finally, for certain dynamical systems with an isolated fixed point, we construct a finite-dimensional Koopman operator that also advances the state directly.  
For the examples presented, this procedure is closely related to the Carleman linearization~\cite{carleman1932application,steeb1980non,kowalski1991nonlinear}, which has extensions to nonlinear control~\cite{Brockett1976automatica,banks1992infinite,svoronos1994discretization}.  
Afterward, it is possible to develop a \emph{nonlinear} Koopman operator optimal control (KOOC) law, even for nonlinear fixed points, using techniques from linear optimal control theory.

\section{Background on Koopman analysis}\label{sec:background}
Consider a continuous-time dynamical system, given by:
\begin{eqnarray}
\frac{d}{dt}\mathbf{x} = \mathbf{f}(\mathbf{x}),\label{Eq:statesyscont}
\end{eqnarray}
where $\mathbf{x}\in\mathbf{M}$ is an $n$-dimensional state on a smooth manifold $\mathbf{M}$.  
The vector field $\mathbf{f}$ is an element of the tangent bundle $\mathbf{T}\mathbf{M}$ of $\mathbf{M}$, such that $\mathbf{f}(\mathbf{x})\in \mathbf{T}_\mathbf{x}\mathbf{M}$.   
Note that in many cases we dispense with manifolds and choose $\mathbf{M}=\mathbb{R}^n$ and $\mathbf{f}$ a Lipschitz continuous function.

For a given time $t$, we may consider the flow map $\mathbf{F}_t:\mathbf{M}\rightarrow\mathbf{M}$, which maps the state $\mathbf{x}(t_0)$ forward time $t$ into the future to $\mathbf{x}(t_0+t)$, according to:
\begin{eqnarray}
\mathbf{F}_t(\mathbf{x}(t_0)) = \mathbf{x}(t_0+t) = \mathbf{x}(t_0) + \int_{t_0}^{t_0+t}\mathbf{f}(\mathbf{x}(\tau))\,d\tau.
\end{eqnarray}
In particular, this induces a discrete-time dynamical system:
\begin{eqnarray}
\mathbf{x}_{k+1} = \mathbf{F}_{t}(\mathbf{x}_k),\label{Eq:statesysdisc}
\end{eqnarray}
where $\mathbf{x}_k = \mathbf{x}(kt)$.   
In general, discrete-time dynamical systems are \emph{more} general than continuous time systems, but we choose to start with continuous time for illustrative purposes.

We also define a real-valued observable function $g:\mathbf{M}\rightarrow\mathbb{R}$, which is an element of an infinite-dimensional Hilbert space.  
Typically, the Hilbert space is given by the Lebesque square-integrable functions on $\mathbf{M}$; other choices of a measure space are also valid.  

The Koopman operator $\mathcal{K}_{t}$ is an infinite-dimensional linear operator that acts on observable functions $g$ as:
\begin{eqnarray}
\mathcal{K}_t g = g\circ \mathbf{F}_t
\end{eqnarray}
where $\circ$ is the composition operator, so that:
\begin{eqnarray}
\mathcal{K}_t g(\mathbf{x}_k) = g(\mathbf{F}_t(\mathbf{x}_k)) = g(\mathbf{x}_{k+1}).
\end{eqnarray}
In other words, the Koopman operator $\mathcal{K}_t$ defines an infinite-dimensional linear dynamical system that advances the observation of the state $g_k=g(\mathbf{x}_k)$ to the next timestep:
 \begin{eqnarray}
 g(\mathbf{x}_{k+1}) = \mathcal{K}_tg(\mathbf{x}_k).\label{Eq:obsvsysdisc}
 \end{eqnarray}
 Note that this is true for \emph{any} observable function $g$ and for any point $\mathbf{x}_k\in\mathbf{M}$.  
 
In the original paper by Koopman, Hamiltonian fluid systems with a positive density were investigated.  
In this case, the Koopman operator $\mathcal{K}_t$ is unitary, and forms a one-parameter family of unitary transformations in Hilbert space.  
The Koopman operator is also known as the composition operator, which is formally the pull-back operator on the space of scalar observable functions~\cite{MarsdenMTAA}.  
The Koopman operator is the dual, or left-adjoint, of the Perron-Frobenius operator, or transfer operator, which is the push-forward operator on the space of probability density functions.   

We may also describe the continuous-time version of the Koopman dynamical system in Eq.~\eqref{Eq:obsvsysdisc} with the infinitesimal generator $\mathcal{K}$ of the one-parameter family of transformations $\mathcal{K}_t$~\cite{MarsdenMTAA} :
\begin{eqnarray}
\frac{d}{dt}g = \mathcal{K}g.\label{Eq:obsvsyscont}
\end{eqnarray}
The linear dynamical systems in Eqs.~\eqref{Eq:obsvsyscont} and \eqref{Eq:obsvsysdisc} are analogous to the dynamical systems in Eqs.~\eqref{Eq:statesyscont} and \eqref{Eq:statesysdisc}, respectively.  
It is important to note that the original state $\mathbf{x}$ may be the observable, and the infinite-dimensional operator $\mathcal{K}_t$ will advance this observable function.  
Again, for Hamiltonian systems, the infinitesimal generator $\mathcal{K}$ is self-adjoint.

\section{Koopman invariant subspaces and exact finite-dimensional models}\label{Sec:SparseKoopman}
As with any vector space, we may choose a basis for Hilbert space and represent our observable function $g$ in this basis.  
For simplicity, let us consider basis observable functions $y_1(\mathbf{x})$, $y_2(\mathbf{x})$, etc., and let a given function $g(\mathbf{x})$ be written in these coordinates as:
\begin{eqnarray}
g = \sum_{k=1}^{\infty}\alpha_ky_k.
\end{eqnarray}
A \emph{Koopman-invariant subspace} is given by $\text{span}\{y_{s_1},y_{s_2},\cdots,y_{s_m}\}$ if all functions $g$ in this subspace, 
\begin{eqnarray}
g=\alpha_1y_{s_1}+\alpha_2y_{s_2}+\cdots+\alpha_my_{s_m},
\end{eqnarray}
remain in this subspace after being acted on by the Koopman operator $\mathcal{K}$:
\begin{eqnarray}
\mathcal{K}g =\beta_1 y_{s_1} + \beta_2 y_{s_2} + \cdots + \beta_m y_{s_m}.
\end{eqnarray}
For functions in these invariant subspaces, it is possible to restrict the Koopman operator to this subspace, yielding a finite-dimensional linear operator $\mathbf{K}$.  $\mathbf{K}$ acts on a vector space $\mathbb{R}^m$, with the coordinates given by the values of $y_{s_k}(\mathbf{x})$.  
This induces a finite-dimensional linear system, as in Eqs.~\eqref{Eq:obsvsysdisc} and~\eqref{Eq:obsvsyscont}.  
Koopman eigenfunctions $\varphi$, such that $\mathcal{K}\varphi = \lambda\varphi$, generate invariant subspaces; however, it may or may not be possible to invert these functions to recover the original state $\mathbf{x}$.

\begin{figure}
\begin{center}
\vspace{-.15in}
\includegraphics[width=.75\textwidth]{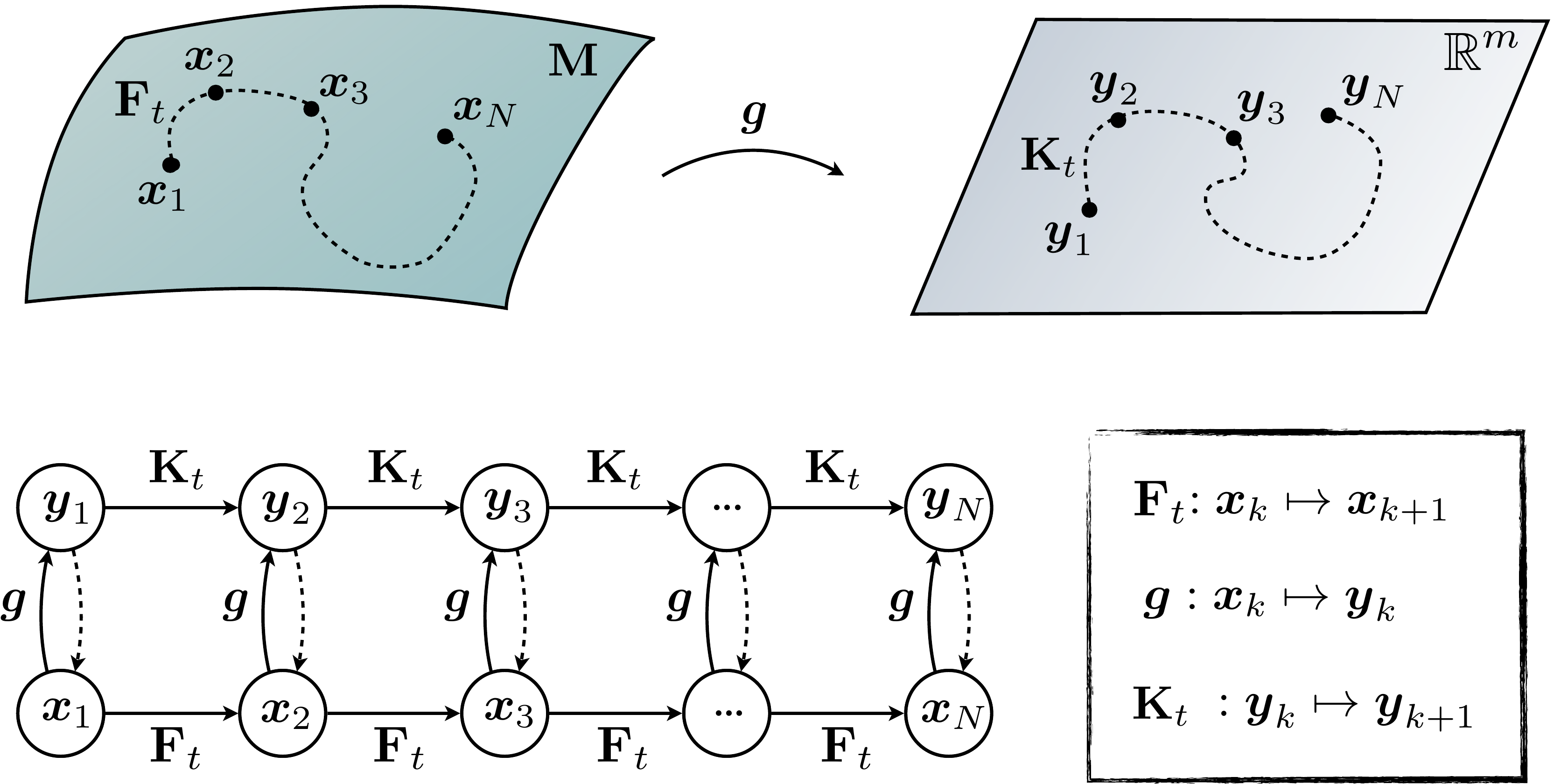}
\end{center}
\vspace{-.15in}
\caption{Schematic illustrating the Koopman operator for nonlinear dynamical systems.  
The dashed lines from $y_k\rightarrow x_k$ indicate that we would like to be able to recover the original state.}\label{fig:overview}
\end{figure}

For control, we may seek Koopman-invariant subspaces that include the original state variables $x_1$, $x_2$, $\cdots$, $x_n$.  
The Koopman operator restricted to this subspace is finite-dimensional, linear, and it advances the original state dynamics, as well as the other observables in the subspace, as shown in Fig.~\ref{fig:overview}. 
These Koopman-invariant subspaces may be identified using data-driven methods, as discussed in Sec.~\ref{sec:SINDY}.
In the following sections, we will show that including the state in our observable subspace is rather restrictive, and it is not possible for the vast majority of nonlinear systems.  
In fact, it is impossible to determine a finite-dimensional Koopman-invariant subspace that includes the original state variables for any system with multiple fixed points or any more general attractors.  
This is because all finite-dimensional linear systems have a single fixed point, and cannot be topologically conjugate to a system with multiple fixed points. 
This does not, however, preclude the identification of Koopman-invariant subspaces spanned by Koopman eigenfunctions $\varphi$, which may provide useful intrinsic coordinates~\cite{Williams2015epl}.  
In fact, it is possible to establish topological conjugacy of the entire basin of attraction of a stable or unstable fixed point or periodic orbit with an associated linear system through the Koopman operator, as shown in~\cite{Lan2013physd,Mezic2013arfm}.  
It may be possible to invert these coordinates to recover the states, although determining eigenfunctions and inverting them to obtain the state may both be challenging. 

For the original state variables $x_1$, $x_2$, $\cdots$, $x_n$ to be included in the Koopman-invariant subspace, then the nonlinear right hand side function $\mathbf{f}$ must also be in this subspace:  
\begin{eqnarray}
\frac{d}{dt}\mathbf{x} = \mathbf{f}(\mathbf{x}) \quad\Longrightarrow\quad \frac{d}{dt}\begin{bmatrix} x_1\\ x_2 \\ \vdots \\ x_n\end{bmatrix} = \begin{bmatrix} f_1(x_1, x_2, \cdots , x_n)\\ f_2(x_1, x_2, \cdots, x_n)\\ \vdots\\ f_n(x_1, x_2, \cdots, x_n)\end{bmatrix}.\label{Eq:dynamicsexpanded}
\end{eqnarray}
The state variables $\mathbf{x}$ form the first $n$ observable functions $y_{s_1}=x_1$, $y_{s_2}=x_2$, $\cdots$, $y_{s_n}=x_n$, and the remaining $m-n$ observables are nonlinear functions required to represent the terms in $\mathbf{f}$.  
If it is possible to represent each term $f_k$ as a combination of observable functions in the subspace,
\begin{eqnarray}
f_k(x_1,x_2, \cdots, x_n) = c_{k,1}y_{s_1} + c_{k,2}y_{s_2} + \cdots + c_{k,m}y_{s_m},\label{Eq:dynamicsasobsv}
\end{eqnarray}
then we may write the first $n$ rows of the Koopman-induced dynamical system as:
\begin{eqnarray}
\frac{d}{dt}\begin{bmatrix} y_1\\ y_2 \\ \vdots\\ y_n \\y_{n+1} \\  \vdots \\ y_m\end{bmatrix} = 
\begin{bmatrix}
c_{1,1} & c_{1,2} & \cdots & c_{1,n} & c_{1,n+1} & \cdots & c_{1,m}\\
c_{2,1} & c_{2,2} & \cdots & c_{2,n} & c_{2,n+1} & \cdots  & c_{2,m}\\
\vdots & \vdots & \ddots & \vdots & \vdots & \ddots & \vdots \\ 
c_{n,1} & c_{n,2} & \cdots & c_{n,n} & c_{n,n+1} & \cdots  & c_{n,m}\\
? & ? & \cdots & ? & ? & \cdots & ?\\
\vdots & \vdots & \ddots & \vdots & \vdots & \ddots & \vdots \\ 
? & ? & \cdots & ? & ? & \cdots & ?
\end{bmatrix}
\begin{bmatrix} y_1\\ y_2 \\ \vdots\\ y_n \\ y_{n+1}\\ \vdots \\ y_m\end{bmatrix}\label{Eq:obsvupdate}.
\end{eqnarray}
In practice, the last $m-n$ rows may be determined analytically, by successively computing $\frac{d}{dt}y_k$ for $k>n$ and representing these derivatives in terms of other subspace observables.  
Knowing the dynamics $\mathbf{f}$ is essential to choose a relevant observable subspace.  
If the dynamics are known, observables may be derived analytically.  
Alternatively, a least-squares regression may be performed using data, as in the extended DMD~\cite{Williams2015jnls}.  

\subsection{Data-driven sparse identification of nonlinear observable functions}\label{sec:SINDY}
It is clear from Eqs.~\eqref{Eq:dynamicsexpanded}--\eqref{Eq:obsvupdate} that the choice of relevant Koopman observable functions is closely related to the form of the nonlinearity in the dynamics.  
In the case that governing equations are unknown, \emph{data-driven} strategies must be employed to determine useful observable functions.  
A recently developed technique allows for the identification of the nonlinear dynamics in Eq.~\eqref{Eq:dynamicsexpanded}, purely from measurements of the system~\cite{Brunton2015arxiv}.  
The so-called sparse identification of nonlinear dynamics (SINDy) algorithm uses sparse regression~\cite{Tibshirani1996lasso} in a nonlinear function space to determine the relevant terms in the dynamics.  
This may be thought of as a generalization of earlier methods that employ symbolic regression (i.e., genetic programming~\cite{koza1999genetic}) to identify dynamics~\cite{Bongard2007pnas,Schmidt2009science}; a similar method has been used to predict catastrophes in dynamical systems~\cite{Wang:2011}.  
Thus, the SINDy algorithm is an equation-free method~\cite{Kevrekidis2003cms} to identify a dynamical system from data.  
This follows a growing trend to exploit sparsity in dynamics~\cite{Ozolicnvs2013pnas,Schaeffer2013pnas,mackey2014compressive} and dynamical systems~\cite{Bai2014aiaa,Proctor2014epj,Brunton2014siads}.

For simplicity in connecting the SINDy algorithm with dynamic mode decomposition (DMD), we consider discrete-time systems as in Eq.~\eqref{Eq:statesysdisc}, although the algorithm applies equally well to continuous-time systems.  
In the SINDy algorithm, measurements of the state $\mathbf{x}$ of a dynamical system are collected, and these measurements are augmented into a larger vector $\bTheta(\bx)$ which contains candidate functions $y_{c_k}(\mathbf{x})$ for the right-hand side dynamics $\mathbf{F}_t(\bx)$ in Eq.~\eqref{Eq:statesysdisc}:
\begin{eqnarray}
\bTheta(\bx) = \begin{bmatrix}y_{c_1}(\bx)\\ y_{c_2}(\bx)\\ \vdots\\ y_{c_m}(\bx)\end{bmatrix}.
\end{eqnarray}
Often, we will choose the first $n$ functions to be the original state variables, $y_{c_k}(\bx) =x_k$, so that the state $\bx$ in $\bTheta(\bx)$.  
Then, we write the following matrix system of equations:
\begin{eqnarray}
\begin{bmatrix}\vline & \vline &  & \vline \\ 
\bx_2  & \bx_3 & \cdots & \bx_M\\
\vline & \vline &  & \vline \end{bmatrix}
 = 
 \begin{bmatrix} \rule[.5ex]{3.em}{0.4pt} & \bxi_1^T & \rule[.5ex]{3.em}{0.4pt} \\
 & \vdots & \\ 
 \rule[.5ex]{3.em}{0.4pt} & \bxi_n^T & \rule[.5ex]{3.em}{0.4pt}\end{bmatrix}
\begin{bmatrix}\vline & \vline &  & \vline \\ 
\bTheta(\bx_1)  & \bTheta(\bx_2) & \cdots & \bTheta(\bx_M)\\
\vline & \vline &  & \vline \end{bmatrix}.
\label{Eq:SINDY}
\end{eqnarray}
This may be written in matrix short-hand as:
\begin{eqnarray}
\bX'= \bXi^T \bTheta(\bX).\label{Eq:SINDY2}
\end{eqnarray}

The functions in $\bTheta(\bx)$ are candidate terms in the right hand side dynamics $\mathbf{F}_t$, and they will also be candidate observable functions.  
The row vectors $\bxi_k^T$ determine which nonlinear terms in $\bTheta(\bx)$ are active in the $k$-th row of $\mathbf{F}_t$; typically, $\bxi_k$ will be a sparse vector, since only a few terms are active in the right hand side of many dynamical systems of interest.  
In this case, we may use sparse regression to solve for each sparse row $\bxi_k^T$.  
Afterward, the sparse matrix $\bXi^T$ yields a nonlinear discrete-time model for Eq.~\eqref{Eq:statesysdisc}, obtained purely from data:
\begin{eqnarray}
\bx_{k+1} = \bXi^T\bTheta(\bx_k).
\label{Eq:SINDYdynamics}
\end{eqnarray}

With the active terms in the nonlinear dynamics identified as the nonzero entries in the rows of $\bXi^T$, it is possible to include these functions in the Koopman subspace.  
Note that in the original SINDy algorithm, the transpose of Eq.~\eqref{Eq:SINDY} was used so that the rows of $\bXi$ become sparse column vectors, establishing a closer resemblance to sparse regression and compressed sensing formulations.  
Again, either discrete-time or continuous time formulations may be used.  
After a reduced observable subspace has been identified, we may re-apply the SINDy Algorithm: 
\begin{eqnarray}
\bTheta_{\text{ref}}(\bX') = \bXi_{\text{aug}}^T \bTheta(\bX)
\end{eqnarray}
where $\bTheta_{\text{ref}}$ is a refined set of candidate observable functions that are active in Eq.~\eqref{Eq:SINDYdynamics}.  
The additional rows of $\bXi_{\text{aug}}$ determine how these observable functions advance as a linear combination of other observable functions.  
This procedure may be iterated until the subspace converges.  
Also, the $\ell_1$ sparse regularization may be omitted in these regressions.  

\subsubsection{Connections to dynamic mode decomposition (DMD)}
In the case that $\bTheta(\bx) = \bx$, the problem in Eq.~\eqref{Eq:SINDY2} reduces to the standard DMD problem:
\begin{eqnarray}
\bX' = \bXi\bX.
\end{eqnarray}
In the standard DMD algorithm, a solution $\bXi$ is obtained that minimizes the sum-square error:
\begin{eqnarray}
\bXi = \argmin_{\tilde{\bXi}} \|\bX'-\tilde{\bXi}\bX\|_F,
\end{eqnarray}
where $\|\cdot\|_F$ is the Frobenius norm.  
This is generally obtained by computing the pseudo-inverse of $\bX$ using the singular value decomposition (SVD).

\section{Systems with Koopman-invariant subspaces containing the state}
Here, we construct a family of nonlinear dynamical systems where it is possible to find a Koopman-invariant subspace that also includes the original state variables as observable functions. 
These systems necessarily only have a single isolated fixed point, as there is no finite-dimensional linear system that can represent multiple fixed points or more general attractors.  
It is, however, possible to obtain linear representations of entire basins of attractions of certain fixed points using eigenfunction coordinates~\cite{Lan2013physd,Williams2015jnls}, where it may be possible to invert to find the state.  
However, these are still not global descriptions, and finding these eigenfunctions and inverting to recover the state remains an open challenge for most systems.
All of the examples below exhibit polynomial nonlinearities that give rise to polynomial slow or fast manifolds. 

\subsection{Continuous-time formulation}
Consider a continuous-time dynamical system with a polynomial slow manifold, given by
\begin{eqnarray}
\frac{d}{dt}\begin{bmatrix} x_1 \\ x_2\end{bmatrix} & = & 
\begin{bmatrix} \mu x_1\\ \lambda\left(x_2-P(x_1)\right)\end{bmatrix},
\end{eqnarray}
where $P(x)$ is a polynomial function.  
If $\lambda \ll |\mu| < 0$, then $x_2=P(x_1)$ is an asymptotically attracting slow manifold.  
This system has a single fixed point at the origin $x_1=x_2=0$.  
We will show that there always exists a finite-dimensional linear system that is given by the closure of the Koopman operator on an observable subspace spanned by the states $x_1$, $x_2$ and the active polynomial terms in $P(x_1)$.

First, consider a single monomial term given by $P(x) = x^N$.  
Thus, we would augment the state with an observable function $x^N$, so that:
\begin{eqnarray}
\mathbf{y} = \begin{bmatrix} y_1 \\ y_2 \\ y_3 \end{bmatrix} = \begin{bmatrix} x_1 \\ x_2 \\ x_1^N\end{bmatrix}.
\end{eqnarray}
Now, the first two terms for $\frac{d}{dt}y_1=\mu y_1$ and $\frac{d}{dt}y_2=\lambda y_2 - \lambda y_3$ are linearly related to the entries of $\mathbf{y}$.  
Finally, to determine $\frac{d}{dt}y_3$, we need only apply the chain rule:
\begin{eqnarray}
\frac{d}{dt}y_3 = \frac{d}{dt}x_1^N = Nx_1^{N-1}\frac{d}{dt}x_1 = \mu Nx_1^N = \mu N y_3.\label{Eq:chain}
\end{eqnarray}
This is closely related to Carleman linearization~\cite{carleman1932application,steeb1980non,kowalski1991nonlinear}. Thus, the system simplifies as:
\begin{eqnarray}
\frac{d}{dt}\begin{bmatrix} y_1\\ y_2 \\ y_3\end{bmatrix} & = & 
\begin{bmatrix} \mu & 0 & 0 \\ 
0 & \lambda & -\lambda \\ 
0 & 0 & \mu N \end{bmatrix}
\begin{bmatrix} y_1 \\ y_2 \\ y_3\end{bmatrix}.
\end{eqnarray}

For more general polynomials, given by $P(x) = a_1x^{N_1} + a_2 x^{N_2} + \cdots + a_M x^{N_M}$, we have:
\begin{eqnarray}
\hspace{-.2in}\begin{bmatrix} y_1 \\ y_2 \\ y_3 \\ y_4 \\ \vdots \\ y_{M+2}\end{bmatrix}= \begin{bmatrix} x_1 \\ x_2 \\ x_1^{N_1} \\ x_1^{N_2} \\ \vdots \\ x_1^{N_M}\end{bmatrix}\quad\Longrightarrow\quad
\frac{d}{dt}\begin{bmatrix} y_1 \\ y_2 \\ y_3 \\ y_4 \\ \vdots \\ y_{M+2}\end{bmatrix}&=&
\begin{bmatrix}
\mu & 0 & 0 & 0 & \cdots & 0 \\ 
0 & \lambda & -a_1\lambda & -a_2\lambda & \cdots & -a_M\lambda\\
0 & 0 & \mu N_1 & 0 & \cdots & 0 \\
0 & 0 & 0 & \mu N_2 & \cdots & 0 \\
\vdots & \vdots & \vdots & \vdots & \ddots & \vdots \\ 
0 & 0 & 0 & 0 & \cdots & \mu N_M
\end{bmatrix}
\begin{bmatrix}y_1 \\ y_2 \\ y_3 \\ y_4 \\ \vdots \\ y_{M+2}\end{bmatrix}.\label{Eq:slowpoly}
\end{eqnarray}
This expression is finite-dimensional and linear, and it advances the original state $\mathbf{x}$ forward exactly, even though the governing dynamics are nonlinear.

\subsubsection{Continous-time examples}
\begin{figure}
\begin{center}
\vspace{-.1in}
\begin{overpic}[width=.75\textwidth]{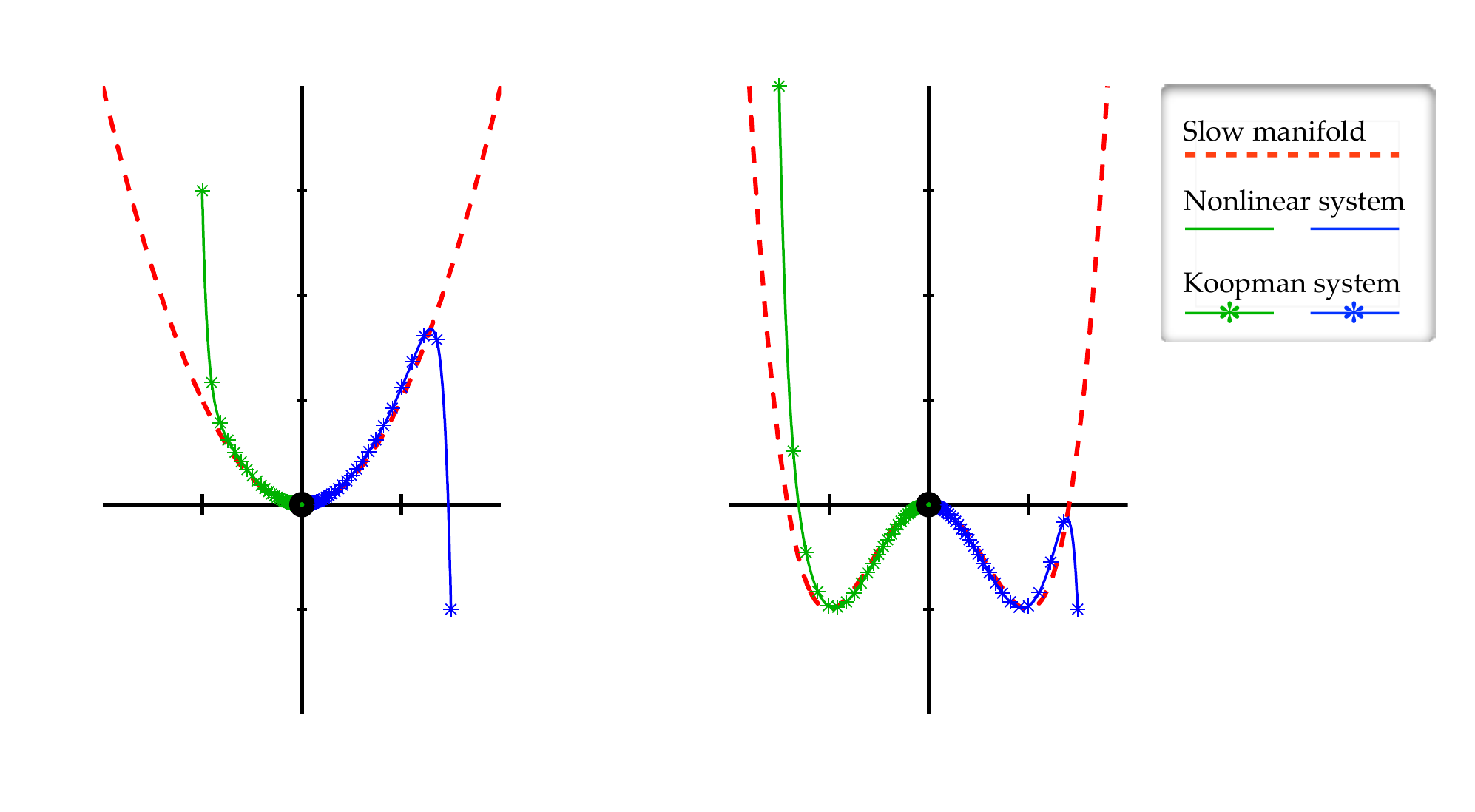}
\put(2,36){(a)}
\put(45,36){(b)}
\end{overpic}
\end{center}
\vspace{-.2in}
\caption{Illustration of two examples with a slow manifold .  In both cases, $\mu=-0.05$ and $\lambda = -1$.}\label{Fig:slowpoly}
\end{figure}

Here, we consider two examples with slow manifolds, which are illustrated in Fig.~\ref{Fig:slowpoly}.  
The first system, with quadratic attracting manifold $x_2=x_1^2$, is given by:
\begin{eqnarray}
\left.\begin{split}
\dot x_1 &= \mu x_1\\
\dot x_2 &= \lambda(x_2-x_1^2)
\end{split}\right\} \quad\Longrightarrow\quad 
\frac{d}{dt}\begin{bmatrix} y_1\\ y_2\\ y_3\end{bmatrix} = 
\begin{bmatrix} \mu & 0 & 0 \\ 0& \lambda& - \lambda\\ 0 & 0 & 2\mu\end{bmatrix}
\begin{bmatrix}y_1\\ y_2\\ y_3\end{bmatrix}
\quad \text{for}\quad 
\begin{bmatrix}y_1\\ y_2\\ y_3\end{bmatrix} = 
\begin{bmatrix}x_1\\ x_2\\x_1^2\end{bmatrix}\label{Eq:QuadraticAttractor}
\end{eqnarray}
and the second system, with quartic attracting manifold $x_2=x_1^4-2x_1^2$, is given by:
\begin{eqnarray}
\left.\begin{split}
\dot x_1 &= \mu x_1\\
\dot x_2 &= \lambda(x_2-x_1^4+2x_1^2)
\end{split}\right\} \quad\Longrightarrow\quad 
\frac{d}{dt}\begin{bmatrix} y_1\\ y_2\\ y_3\\ y_4\end{bmatrix} = 
\begin{bmatrix} \mu & 0 & 0  & 0 \\ 0& \lambda& 2\lambda & -\lambda\\ 0 & 0 & 2\mu & 0 \\ 0 & 0 & 0 & 4\mu\end{bmatrix}
\begin{bmatrix}y_1\\ y_2\\ y_3\\ y_4\end{bmatrix}
\quad \text{for}\quad 
\begin{bmatrix}y_1\\ y_2\\ y_3\\ y_4\end{bmatrix} = 
\begin{bmatrix}x_1\\ x_2\\x_1^2\\ x_1^4\end{bmatrix}.
\end{eqnarray}

\begin{figure}
\begin{center}
\begin{minipage}{.6\textwidth}
\hspace{-.2in}
\begin{overpic}[width=\textwidth]{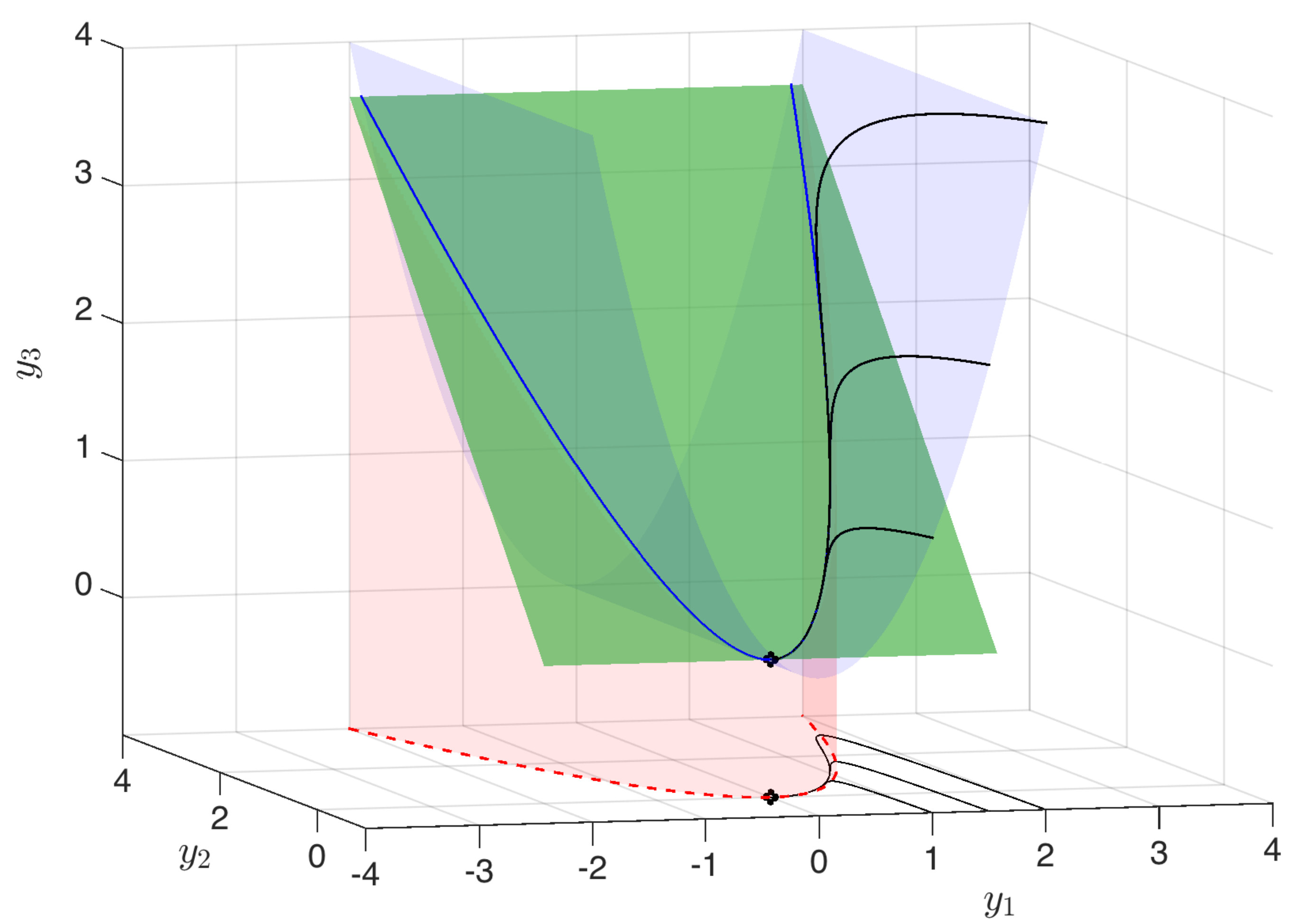}
\end{overpic}
\end{minipage}
\begin{minipage}{.35\textwidth}
\caption{Visualization of three-dimensional linear Koopman system from Eq.~\eqref{Eq:QuadraticAttractor} along with projection of dynamics onto the $x_1$-$x_2$ plane.  The attracting slow manifold is shown in red, the constraint $y_3=y_1^2$ is shown in blue, and the slow unstable subspace of Eq.~\eqref{Eq:QuadraticAttractor} is shown in green.  Black trajectories of the linear Koopman system in $\mathbf{y}$ project onto trajectories of the full nonlinear system in $\mathbf{x}$ in the $y_1$-$y_2$ plane.  Here, $\mu=-0.05$ and $\lambda=1$.  Figure is reproduced with Code~\ref{Code:quad_manifold}.}\label{Fig:quad_manifold}
\end{minipage}
\vspace{-.25in}
\end{center}
\end{figure}

To understand the embedding of a nonlinear dynamical system in a higher-dimensional observable subspace, in which the dynamics are linear, consider the system with quadratic attracting manifold from Eq.~\eqref{Eq:QuadraticAttractor}.  
The full three-dimensional Koopman observable vector space is visualized in Fig.~\ref{Fig:quad_manifold}.  
Trajectories that start on the invariant manifold $y_3=y_1^2$, visualized by the blue surface, are constrained to stay on this manifold.  
There is a \emph{slow} subspace, spanned by the eigenvectors corresponding to the slow eigenvalues $\mu$ and $2\mu$; this subspace is visualized by the green surface.  
Finally, there is the original asymptotically attracting manifold of the original system, $y_2=y_1^2$, which is visualized as the red surface.  
The blue and red parabolic surfaces always intersect in a parabola that is inclined at a $45^{\circ}$ angle in the $y_2$-$y_3$ direction.  
The green surface approaches this $45^{\circ}$ inclination as the ratio of fast to slow dynamics become increasingly large.  
In the full three-dimensional Koopman observable space, the dynamics are given by a stable node, with trajectories rapidly attracting onto the green subspace and then slowly approaching the fixed point.  

\subsubsection{Intrinsic coordinates defined by eigen-observables of the Koopman operator}
The left eigenvectors of the Koopman operator yield Koopman eigenfunctions (i.e., eigenobservables).  
The Koopman eigenfunctions of Eq.~\eqref{Eq:QuadraticAttractor} corresponding to eigenvalues $\mu$ and $\lambda$ are:
\begin{eqnarray}
\varphi_{\mu}=x_1, \quad\text{and}\quad\varphi_{\lambda} = x_2-bx_1^2 \quad\text{with}\quad b =\frac{\lambda}{\lambda-2\mu}.
\end{eqnarray}
The constant $b$ in $\varphi_\lambda$ captures the fact that for a finite ratio $\lambda/\mu$, the dynamics only shadow the asymptotically attracting slow manifold $x_2=x_1^2$, but in fact follow neighboring parabolic trajectories.  
This is illustrated more clearly by the various surfaces in Fig.~\ref{Fig:quad_manifold} for different ratios $\lambda/\mu$.  

In this way, a set of intrinsic coordinates may be determined from the observable functions defined by the left eigenvectors of the Koopman operator on an invariant subspace.  
Explicitly, 
\begin{eqnarray}
\varphi_{\alpha}(\bx) = \bxi_{\alpha}\by(\bx), \quad\text{where}\quad \bxi_{\alpha}\mathbf{K} = \alpha\bxi_{\alpha}.
\end{eqnarray}
These eigen-observables define observable subspaces that remain invariant under the Koopman operator, even after coordinate transformations.  
As such, they may be regarded as intrinsic coordinates~\cite{Williams2015epl} on the Koopman-invariant subspace.  
As an example, consider the system from Eq.~\eqref{Eq:QuadraticAttractor}, but written in a coordinate system that is rotated by $45^{\circ}$:
\begin{eqnarray}
\left.\begin{split}
\eta & = x + y\\ 
\xi  & = x - y
\end{split}
\quad\text{ and } \quad
\begin{split}
x & =\left(\eta+\xi\right)/2\\ 
y  & =\left(\eta-\xi\right)/2
\end{split}\right\} \quad \Longrightarrow\quad
 \begin{split}
 \frac{d}{dt}\eta & = \frac{\mu}{2}(\eta+\xi) +\frac{\lambda}{2}(\eta-\xi)-\frac{\lambda}{4}(\eta+\xi)^2\\
 \frac{d}{dt}\xi & = \frac{\mu}{2}(\eta+\xi) -\frac{\lambda}{2}(\eta-\xi)+\frac{\lambda}{4}(\eta+\xi)^2\\
 \end{split}
\end{eqnarray}
The original eigenfunctions, written in the new coordinate systems are:
\begin{eqnarray*}
&&\varphi_\mu(\eta,\xi) = \frac{\eta+\xi}{2}\\
&&\varphi_{\lambda}(\eta,\xi) = \frac{\eta-\xi}{2} - \frac{\lambda}{\lambda-2\mu}\frac{(\eta+\xi)^2}{4}.
\end{eqnarray*}
It is easy to verify that these remain eigenfunctions:
\begin{eqnarray*}
\frac{d}{dt} \varphi_{\mu} &=& \frac{\dot\eta+\dot\xi}{2} = \mu \frac{\eta+\xi}{2} = \mu\varphi_{\mu}\\
\frac{d}{dt}\varphi_{\lambda} &=& \frac{\dot\eta-\dot\xi}{2} - \frac{\lambda}{\lambda-2\mu}\frac{2(\eta+\xi)(\dot\eta+\dot\xi)}{4}\\
&=&\lambda\left[\frac{\eta-\xi}{2} - \frac{\lambda}{\lambda-2\mu}\frac{(\eta+\xi)^2}{4}\right] = \lambda\varphi_{\lambda}.
\end{eqnarray*}
In fact, in this new coordinate system, it is possible to write the Koopman subspace system:
\begin{eqnarray}
\frac{d}{dt}\begin{bmatrix}\eta\\ \xi\\ \varphi_{\lambda}\end{bmatrix} = 
\begin{bmatrix}
\frac{3\mu}{2} & -\frac{\mu}{2} & (\lambda-2\mu)\\
-\frac{\mu}{2} & \frac{3\mu}{2} & -(\lambda-2\mu)\\
0 & 0 & \lambda
\end{bmatrix}
\begin{bmatrix}\eta\\ \xi\\ \varphi_{\lambda}\end{bmatrix} .
\end{eqnarray}

\subsection{Discrete-time formulation}
A related formulation for discrete-time systems is given by:
\begin{eqnarray}
\begin{bmatrix}x_1 \\ x_2\end{bmatrix}_{k+1} = \begin{bmatrix}\mu & 0 \\ 0 & \lambda\end{bmatrix}\begin{bmatrix}x_1 \\ x_2 \end{bmatrix}_k + \begin{bmatrix} 0 \\ (1-\lambda)P([x_1]_k)\end{bmatrix}.
\end{eqnarray}
This system will also converge asymptotically to a slow manifold given by $x_2=P(x_1)$ when $|\lambda|\ll |\mu |$ and $|\lambda|<1$.  
A similar argument can be made to that given in Eq.~\eqref{Eq:chain} and Eq.~\eqref{Eq:slowpoly}, but with $\mu^N$ replacing $N\mu$, since:
\begin{eqnarray}
[x_1^N]_{k+1} = \left([x_1]_{k+1}\right)^N = \left(\mu [x_1]_k\right)^N = \mu^N [x_1^N]_k.
\end{eqnarray}
Thus, for discrete-time systems, the update is given by:
\begin{eqnarray}\nonumber
\begin{bmatrix} y_1 \\ y_2 \\ y_3 \\ y_4 \\ \vdots \\ y_{M+2}\end{bmatrix}= \begin{bmatrix} x_1 \\ x_2 \\ x_1^{N_1} \\ x_1^{N_2} \\ \vdots \\ x_1^{N_M}\end{bmatrix}~~\Longrightarrow~~
\begin{bmatrix} y_1 \\ y_2 \\ y_3 \\ y_4 \\ \vdots \\ y_{M+2}\end{bmatrix}_{k+1}=
\begin{bmatrix}
\mu & 0 & 0 & 0 & \cdots & 0 \\ 
0 & \lambda & a_1(1-\lambda) & a_2(1-\lambda) & \cdots & a_M(1-\lambda)\\
0 & 0 & \mu^N_1 & 0 & \cdots & 0 \\
0 & 0 & 0 & \mu^N_2 & \cdots & 0 \\
\vdots & \vdots & \vdots & \vdots & \ddots & \vdots \\ 
0 & 0 & 0 & 0 & \cdots & \mu^N_M
\end{bmatrix}
\begin{bmatrix}y_1 \\ y_2 \\ y_3 \\ y_4 \\ \vdots \\ y_{M+2}\end{bmatrix}_k.\label{Eq:discreteslowpoly}
\end{eqnarray}

\subsubsection{Discrete-time example}\label{Sec:RowleyWilliams}
The case of a polynomial slow manifold is inspired by an example from Tu et al.~\cite{Tu2014jcd}:
\begin{eqnarray}
\begin{bmatrix} x_1\\ x_2\end{bmatrix} \mapsto \begin{bmatrix} \lambda x_1 \\ \mu x_2 + (\lambda^2-\mu)x_1^2\end{bmatrix}.\label{Eq:Rowley}
\end{eqnarray}
In this case, there is a polynomial stable manifold $x_2=x_1^2$.   
Thus, they suggest the following observable variables, which are intrinsic coordinates for the dynamics:
\begin{eqnarray}
\begin{bmatrix}y_1\\ y_2\end{bmatrix} = \begin{bmatrix}x_1\\ x_2-x_1^2\end{bmatrix} \quad\Longrightarrow\quad \begin{bmatrix}y_1\\ y_2\end{bmatrix}_{k+1} = \begin{bmatrix}\lambda & 0 \\ 0 & \mu \end{bmatrix}\begin{bmatrix}y_1\\ y_2\end{bmatrix}_{k}.
\end{eqnarray}

In our framework above, if the correct intrinsic variables were unknown, they could be discovered by writing the system as:
\begin{eqnarray}
\begin{bmatrix} y_1 \\ y_2 \\ y_3\end{bmatrix} = \begin{bmatrix} x_1 \\ x_2 \\ x_1^2\end{bmatrix}\quad\Longrightarrow\quad 
\begin{bmatrix} y_1\\ y_2\\ y_3\end{bmatrix}_{k+1} = \begin{bmatrix} \lambda & 0 & 0 \\ 0 & \mu & (\lambda^2-\mu)\\ 0 & 0 & \lambda^2\end{bmatrix}\begin{bmatrix} y_1\\ y_2\\ y_3\end{bmatrix}_k.
\end{eqnarray}
Finally, in this observable function coordinate system, the left eigenvectors are:
\begin{eqnarray}
{\xi_1} = \begin{bmatrix} 1 \\ 0 \\ 0\end{bmatrix}~\Longrightarrow~ \varphi_1(\mathbf{x})=x_1, ~~~~ \xi_2=\begin{bmatrix} 0\\ 0\\ 1\end{bmatrix}~\Longrightarrow~ \varphi_2(\mathbf{x})=x_1^2,~~~~ \xi_3 = \begin{bmatrix} 0 \\ 1\\ -1\end{bmatrix}~\Longrightarrow~ \varphi_3(\mathbf{x})=x_2-x_1^2
\end{eqnarray}
corresponding to the eigenvalues $\lambda_1=\lambda$, $\lambda_2=\lambda^2$ and $\lambda_3=\mu$.  
These eigenvectors diagonalize the system and define the intrinsic coordinates.  

\section{Koopman operator optimal control}
A long held hope of Koopman operator theory is that it would provide insights into the control of nonlinear systems.  
Here, we present results of designing control laws using linear control theory on the truncated Koopman operator; these Koopman operator optimal controllers (KOOCs) then induce a nonlinear controller on the state-space that dramatically outperforms optimal control on the linearized fixed point.  

This is only a brief introduction to the theory of Koopman optimal control, and there are numerous extensions that must be developed and explored.  
There are existing connections between DMD and control systems~\cite{Proctor2014arxiv}, and there are ongoing efforts to extend this to the Koopman operator framework.  
There are a number of systems where it is not clear how to use the Koopman linear operator for control, and these will be briefly outlined below.  
In addition, there are alternative nonlinear control methods related to Carleman linearization~\cite{Brockett1976automatica,banks1992infinite,svoronos1994discretization} that may be connected to Koopman operator control.  
Moreover, we have not yet proven the nonlinear optimality of these new controllers, but the numerical performance is striking.  

\subsection{Simple motivating example}
As a motivating example, consider the nonlinear system in Eq.~\eqref{Eq:QuadraticAttractor}, but with the stability of the $x_2$ direction reversed (i.e., $\lambda=1$ instead of $\lambda=-1$), and modified to include actuation on the second state:
\begin{eqnarray}
\frac{d}{dt}\begin{bmatrix}x_1\\ x_2\end{bmatrix} & = &
 \begin{bmatrix}\mu & 0 \\ 0 & \lambda\end{bmatrix}\begin{bmatrix}x_1\\ x_2\end{bmatrix} + 
 \begin{bmatrix} 0 \\  -\lambda x_1^2\end{bmatrix} + \begin{bmatrix} 0\\ 1\end{bmatrix}u,\label{Eq:controller1}
\end{eqnarray}
with $\mu=-.1$ and $\lambda=1$.  
Again, this may be put into a Koopman formalism as:
\begin{eqnarray}
\frac{d}{dt}\begin{bmatrix} y_1\\ y_2\\ y_3\end{bmatrix} & = & 
\begin{bmatrix}
\mu & 0 & 0 \\
0 & \lambda & -\lambda\\
0 & 0 & 2\mu\end{bmatrix}\begin{bmatrix} y_1\\ y_2\\ y_3\end{bmatrix} + 
\begin{bmatrix} 0\\ 1\\ 0\end{bmatrix} u.\label{Eq:controller2}
\end{eqnarray}

Now, let us assume that we have a quadratic cost function, as in the linear-quadratic-regulator (LQR) control framework:
\begin{eqnarray}
J = \int_0^\infty \mathbf{x}^T(\tau)\mathbf{Q}\mathbf{x}(\tau) + \mathbf{u}(\tau)^T\mathbf{R}\mathbf{u}(\tau)\,d\tau,
\end{eqnarray}
where $\mathbf{Q}$ weighs the cost of deviations of the state $\mathbf{x}$ from the origin and $\mathbf{R}$ weighs the cost of control expenditure.  
For now, we will consider the following $\mathbf{Q}$ and $\mathbf{R}$ for simplicity: 
\begin{eqnarray}
\mathbf{Q} = \begin{bmatrix} 1 & 0 \\ 0 & 1\end{bmatrix} \quad \quad\quad\quad R = 1.
\end{eqnarray}
In this way, all state deviations and control expenditures are weighed equally.  

For linear systems, such as the linearization of Eq.~\eqref{Eq:controller1}, it is possible to derive the matrix $\mathbf{C}$ that results in the optimal control law $\mathbf{u}=-\mathbf{C}\mathbf{x}$; this control law is optimal in the sense that it achieves the minimal attainable cost function $J$.  
However, this controller will only be optimal for a small vicinity of the fixed point where linearization is valid.  
Outside this vicinity, when nonlinear terms become large, all guarantees of optimality are lost.  

\begin{figure}
\begin{center}
\begin{overpic}[width=\textwidth]{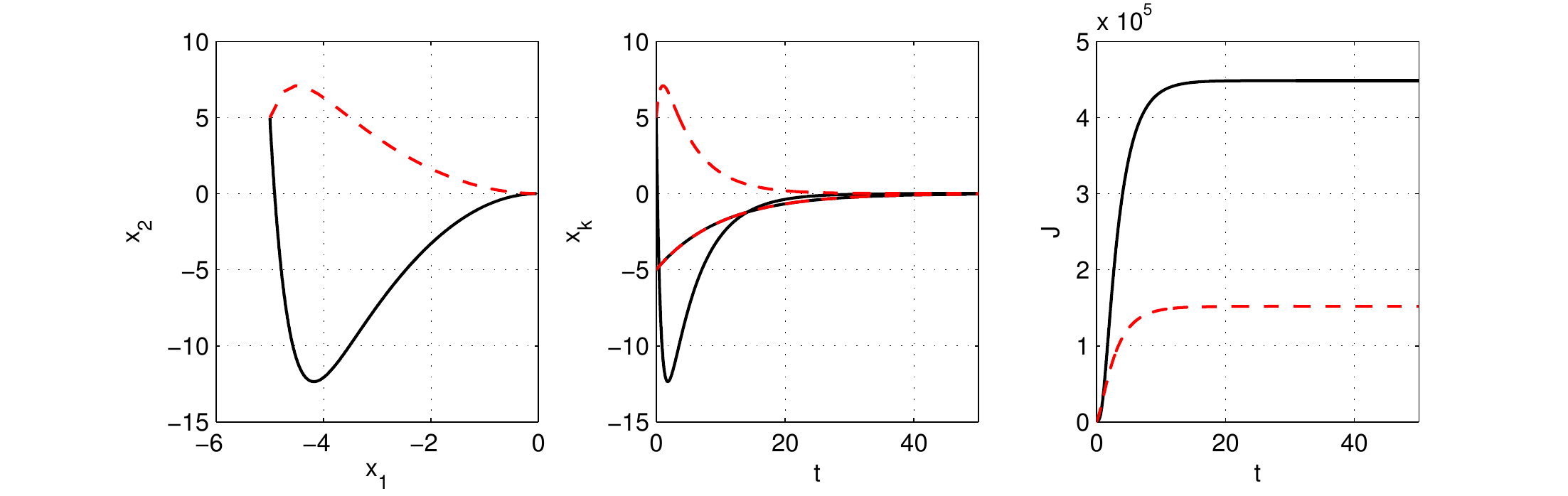}
\end{overpic}
\end{center}
\vspace{-.2in}
\caption{Illustration of LQR control around a nonlinear fixed point using standard linearization (black) and truncated Koopman (red).  The Koopman optimal controller achieves a much smaller overall cost, $J$, approximately $1/3$ of the cost of the standard LQR solution.}\label{Fig:lqr_nonlinear}
\end{figure}

Instead of linearizing near the fixed point and computing the optimal LQR controller, here we use the Koopman linear system in Eq.~\eqref{Eq:controller2}.  
We still have the same cost on the state $\mathbf{x}$, so we use a modified weight matrix $\tilde{\mathbf{Q}}$ given by $\tilde{\mathbf{Q}} = \begin{bmatrix} \mathbf{Q} & 0\\ 0 & 0\end{bmatrix}$ and $\tilde{R}=R$.  
In this way, we may develop an optimal \emph{linear} controller for the Koopman representation of our nonlinear system.  
In this case, the Koopman linear control law, given by $u=\tilde{\mathbf{C}}\mathbf{y}$, may be interpreted as a nonlinear control law on the original state $\mathbf{x}$:
\begin{eqnarray}
u = -\begin{bmatrix}\tilde{K}_1 & \tilde{K}_2\end{bmatrix}\begin{bmatrix} x_1\\ x_2\end{bmatrix} - \tilde{K}_3x_1^2.
\end{eqnarray}
The results of the standard LQR compared with this Koopman operator optimal controller are shown in Fig.~\ref{Fig:lqr_nonlinear}, and the Matlab code is provided in Code~\ref{Code:KOOC}.  
In this example, the KOOC achieves a cost of approximately 1/3 the cost of standard LQR.

\subsection{Limitations of Koopman operator optimal control}
In the current framework, there are a number of limitations to the approach advocated above.  
We will illustrate this on a simple variation on the example above, in which $\mu$ is unstable instead of $\lambda$ and the control input effects the first state $x_1$ instead of $x_2$:
\begin{eqnarray}
\frac{d}{dt}\begin{bmatrix} x_1\\ x_2\end{bmatrix} \begin{bmatrix}\mu & 0 \\ 0 & \lambda\end{bmatrix} \begin{bmatrix}x_1\\ x_2\end{bmatrix} + \begin{bmatrix}0 \\ -\lambda x_1^2\end{bmatrix} + \begin{bmatrix}1\\ 0\end{bmatrix}u,
\end{eqnarray}
with $\mu=.1$ and $\lambda=-1$.  
In this example, it is necessary to move the actuation to the first state $x_1$, otherwise this state will be unstable and uncontrollable.  
What is more troubling, is that the subspace spanned by $x_1$, $x_2$, and $x_1^2$ is no longer Koopman-invariant, since the expression for the time derivative of $y_3=x_1^2$ is more complicated now:
\begin{eqnarray}
\frac{d}{dt} y_3 = 2x_1 \frac{d}{dt}x_1 = 2x_1\left(\mu x_1 + u\right).
\end{eqnarray}
Thus, there is a troublesome extra nonlinear term $x_1u$ in the expression for $\frac{d}{dt}y_3$.  
However, this may not be too large of a problem, considering that we don't weight excursions of $y_3$ in the cost function.  
What is a larger problem, is that the state $y_3$ has a positive eigenvalue $2\mu$, which is uncontrollable.  
Many off-the-shelf packages, such as Matlab, will fail to return an LQR controller for such uncontrollable unstable systems.

\section{Discussion}

In this paper, we have investigated a special choice of Koopman observable functions that form a finite-dimensional subspace of Hilbert space that contains the state in its span and remains invariant under the Koopman operator.  
Any finite collection of Koopman eigenfunctions (i.e., eigen-observables) forms such a Koopman-invariant subspace.  
These Koopman eigenfunctions may be extremely useful, providing intrinsic coordinates for a given nonlinear dynamical system.  
In addition, given such a Koopman-invariant subspace, the Koopman operator restricted to this subspace yields a finite-dimensional linear dynamical system to evolve these observables forward in time.  
However, it is not always clear how to identify relevant Koopman eigenfunctions, either from data or governing equations, or how to invert these coordinates to obtain information about the progression of the underlying state variables.  
Moreover, in many cases with control, the control objectives are defined directly on the state; this is the case in linear quadratic regulator (LQR) control, for example.  
Thus, there is still interest in defining a Koopman invariant subspace that includes the original state variables as observable functions. 

We demonstrate that for a large class of nonlinear systems with a single isolated fixed point, it is possible to obtain such a Koopman-invariant subspace that includes the original state variables.  
Moreover, we present a data-driven technique to identify the relevant Koopman observable functions, leveraging a recent technique that identifies nonlinear dynamical systems in a nonlinear function space using sparse regression; this algorithm is known as the sparse identification of nonlinear dynamics (SINDy).  
We show that the eigen-observables that define this Koopman-invariant subspace may be solved for as left-eigenvectors of the Koopman operator restricted to the subspace in the chosen coordinate system.  
Finally, we demonstrate that the finite-dimensional linear Koopman operator defined on this Koopman-invariant subspace may be used to develop Koopman operator optimal control (KOOC) laws using techniques from linear control theory.  
In particular, we develop an LQR controller using the Koopman linear system, but retaining the cost function defined on the original state.
The resulting control law may be thought of as inducing a nonlinear control law on the state variable, and it dramatically outperforms standard LQR computed on a linearization, reducing the cost expended by a factor of three.
This is extremely promising and may result in significantly improved control laws for systems with normal form expansions near fixed points~\cite{guckenheimer_holmes}.  These expansions are commonly used in astrophysical problems to compute orbits around fixed points~\cite{koon2008dynamical}; for example, the James Webb Space Telescope will orbit the Sun-Earth $L_2$ Lagrange point~\cite{gardner2006james}.  

As is often the case with interesting problems in mathematics, a deeper understanding of one problem opens up a host of other open questions.   
For example, a complete classification of nonlinear systems which admit Koopman-invariant subspaces that include the state variables as observables remains an open and interesting problem.  
It is, however, clear that no system with multiple fixed points, or any periodic orbits or more complex attractors can admit such a finite-dimensional Koopman-invariant subspace containing the state variables explicitly as observables.  
In these cases, another open problem is how to choose observable coordinates so that a finite-rank truncation of the linear Koopman dynamics yields useful results, not just for reconstruction of existing data, but for future state prediction and control.  
Finally, more effort must go into understanding whether or not Koopman operator optimal control laws are optimal in the sense that they minimize the cost function across all possible nonlinear control laws.  

Much of the interest surrounding Koopman analysis and DMD has been centered around the promise of obtaining finite-dimensional linear expressions for nonlinear dynamics.  
In fact, any set of Koopman eigenfunctions span an invariant subspace, where it is possible to obtain an exact and closed finite-dimensional truncation, although finding these nonlinear Koopman eigen-observable functions is challenging.  
Moreover, Koopman invariant subspaces may or may not provide enough information to propagate the underlying state, which is useful for evaluating cost functions in optimal control laws. 
Koopman eigenfunctions provide a wealth of information about the original system, including a characterization of invariant sets such as stable and unstable manifolds, and these may not have simple closed-form representations, but may instead need to be approximated from data.  
There are methods that identify almost invariant sets and coherent structures~\cite{Froyland2005physd,Froyland2009pd} using set oriented methods~\cite{Dellnitz2001book}.  
Related Ulam-Galerkin methods have been used to approximate eigenvalues and eigenfunctions of the Perron-Frobenius operator~\cite{Froyland2014siads}.

To address these challenges, finite-dimensional linear approximations of the Koopman operator from data have been widely explored, and they are valuable in many instances, especially for extracting dynamics on modal coherent structures.  
However, we have shown that it is quite rare for a dynamical system to admit a finite-dimensional Koopman-invariant subspace that includes the state variables explicitly, so that exact linear models to propagate the state dynamics exist only for systems with a single isolated fixed point.  
This implies that approximate truncation of linear Koopman models for nonlinear phenomena with multiple fixed points or more general attractors should be used with care for future-state prediction, especially for off-attractor transients, as well as for the design of control laws.  
There is no free lunch with Koopman analysis of nonlinear systems, as we trade finite-dimensional nonlinear dynamics for infinite-dimensional linear dynamics, with an entirely new host of challenges.

\vspace{-.05in}
\section*{Acknowledgements}
\vspace{-.05in}
The authors would like to thank Clancy Rowley for insightful and constructive comments on an early draft of this paper, especially his observation that the state can often be recovered from the intrinsic eigenfunction coordinates.  
We also thank Igor Mezic for making us aware of Carleman linearization and the extremely interesting literature on operator theoretic control of nonlinear systems.  
SLB acknowledges support from the University of Washington department of Mechanical Engineering and as a Data Science Fellow in the eScience Institute.  BWB acknowledges support from the University of Washington department of Biology and as a Data Science Fellow in the eScience Institute.  JLP thanks Bill and Melinda Gates for their active support of the Institute of Disease Modeling and their sponsorship through the Global Good Fund.  JNK acknowledges support from the U.S. Air Force Office of Scientific Research (FA9550-09-0174). 

\vspace{-.05in}

\small
\bibliographystyle{unsrt}
\bibliography{references}

\section*{Appendix: Systems without Koopman-invariant subspaces that explicitly span the state}
For any system with multiple fixed points, periodic orbits, or atrracting/repelling structures, there is no finite-dimensional Koopman invariant subspace that explicitly includes the state. 
This follows from the fact that these systems cannot be topologically conjugate to a finite-dimensional linear system with a single fixed point.  
It may, however, be possible to obtain a linearization that is valid in an entire basin of attraction of a single fixed point or periodic orbit of a complex system~\cite{Lan2013physd,Williams2015jnls}. 
Moreover, it may be possible to determine an invariant subspace spanned by Koopman eigenfunctions such that it is possible to invert these eigenfunctions to recover the states.  
Incorporating eigenfunction-based Koopman optimal control is an important avenue of future research, as it will open up Koopman optimal nonlinear control to a wider class of important problems.  

\subsection{Example: Logistic map}
Consider the logistic map, given by:
\begin{eqnarray}
x_{k+1} = rx_k(1-x_k).
\end{eqnarray}
Naturally, the observable subspace must include $x$ and $x^2$:
\begin{eqnarray}
\mathbf{y}_{k} = \begin{bmatrix} x\\ x\end{bmatrix}_k \triangleq \begin{bmatrix}x_k\\ x_k^2\end{bmatrix}.
\end{eqnarray}
Writing out the Koopman operator, the first row equation is simple:
\begin{eqnarray}
\mathbf{y}_{k+1} = \begin{bmatrix}x\\ x^2\end{bmatrix}_{k+1} =  \begin{lbmatrix}{2}r & -r \\ \hdashline ? & ?\end{lbmatrix}\begin{bmatrix}x\\ x^2\end{bmatrix}_k,
\end{eqnarray}
but the second row is not obvious.  
To find this expression, expand $(x_{k+1})^2$:
\begin{eqnarray}
x_{k+1}^2  =  \left(rx_k(1-x_k)\right)^2 =  r^2 \left(x_k^2 - 2 x_k^3 + x_k^4\right).
\end{eqnarray}
Thus, we also need cubic and quartic polynomial terms to advance $x^2$.  
Similarly, these terms need polynomials up to sixth and eighth order, respectively, and so on, ad infinitum:
\begin{eqnarray}
&&~~~ x \quad~x^2 \quad~~x^3 \quad\quad x^4 \quad\quad~ x^5 \quad\quad x^6 \quad\quad~x^7 \quad\quad~ x^8 \quad\quad x^9 \quad~~ x^{10}\nonumber\\
\begin{bmatrix} x \\ x^2 \\ x^3 \\ x^4 \\ x^5 \\ \vdots\end{bmatrix}_{k+1}
& = & 
\begin{bmatrix}
r & -r & 0 & 0 & 0 & 0 & 0 & 0 & 0 & 0 & \cdots \\ 
0 & r^2  & -2r^2 & r^2 & 0 & 0 & 0 & 0 & 0 & 0 & \cdots\\
0 & 0 & r^3 & -3 r^3 & 3 r^3 & r^3 & 0 & 0 & 0 & 0 & \cdots\\
0 & 0 & 0 & r^4 & -4r^4 & 6r^4 & -4r^4 & r^4 & 0 & 0 & \cdots\\
0 & 0 & 0 & 0 & r^5 & -5r^5 & 10r^5 & -10r^5 & 5r^5 & -r^5 & \cdots\\
\vdots & \vdots & \vdots & \vdots & \vdots & \vdots & \vdots & \vdots & \vdots & \vdots & \ddots
\end{bmatrix}
\begin{bmatrix} x \\ x^2 \\ x^3 \\ x^4 \\ x^5 \\ \vdots \end{bmatrix}_k.
\end{eqnarray}
It is interesting to note that the rows of this equation are related to the rows of Pascal's triangle, with the $n$-th row scaled by $r^n$, and with the omission of the first row:
\begin{eqnarray}
\begin{bmatrix}x^0\end{bmatrix}_{k+1} = \begin{bmatrix}r^0\end{bmatrix}\begin{bmatrix}x^0\end{bmatrix}_k.
\end{eqnarray}

The representation of the Koopman operator in a polynomial basis is somewhat troubling.  
Not only is there no closure, but the determinant of any finite-rank truncation is very large for $r>1$.  
This illustrates a pitfall associated with naive representation of the infinite dimensional Koopman operator for a simple chaotic system.  
Truncating the system, or performing a least squares fit on an augmented observable vector (i.e., DMD on a nonlinear measurement) yields poor results, with the truncated system only agreeing with the true dynamics for a small handful of iterations.

\subsection{Example: Nonlinear fixed point with a center manifold}
Consider the simple nonlinear system with a single isolated fixed point at the origin:
\begin{eqnarray}
\frac{d}{dt}x=x^2.
\end{eqnarray}
The Carleman linearization approach above would suggest that we augment the observable subspace with the quadratic polynomial $y_2=x^2$, so that:
\begin{eqnarray}
\begin{bmatrix} y_1\\ y_2\end{bmatrix} = \begin{bmatrix} x\\ x^2\end{bmatrix}.
\end{eqnarray}
However, the expression for the time-derivative of $y_2$ requires higher polynomials in $x$:
\begin{eqnarray}
\frac{d}{dt} y_2 = 2x\dot{x} = 2x^3.
\end{eqnarray}
Similarly, if we introduce $y_3=x^3$, then
%
$\frac{d}{dt}y_3 = 3x^2\dot{x}=3x^4$,
%
and so on.  
This results in an infinite Koopman expansion:
\begin{eqnarray}
\frac{d}{dt}\begin{bmatrix}y_1\\ y_2\\ y_3\\ y_4 \\ \vdots\end{bmatrix} & = & 
\begin{bmatrix} 0 & 1 & 0 & 0 &  \cdots \\ 
0 & 0 & 2 & 0  & \cdots\\
0 & 0 & 0 & 3 & \cdots\\
0 & 0 & 0 & 0  & \cdots\\
\vdots & \vdots & \vdots & \vdots  & \ddots\end{bmatrix}
\begin{bmatrix}y_1\\ y_2\\ y_3\\ y_4 \\ \vdots\end{bmatrix} \quad\text{where}\quad\begin{bmatrix}
y_1\\ y_2\\ y_3\\ y_4 \\ \vdots\end{bmatrix} = \begin{bmatrix} x\\ x^2\\ x^3\\ x^4 \\  \vdots\end{bmatrix}.\label{Eq:CarlemanQuad}
\end{eqnarray}
Note that the determinant of any finite-rank truncation of the Koopman operator is $0$, even though the system has finite-time blow up!  For this problem, it is possible to use eigenfunction coordinates to obtain a linear model in terms of an eigenfunction that may be inverted to recover the state\footnote{From a personal communication with C. W. Rowley.}:
\begin{eqnarray}
\varphi(x) = e^{-1/x}\quad\Longrightarrow\quad \frac{d}{dt}\varphi(x) = x^{-2}e^{-1/x}\dot{x} = \varphi(x).\label{Eq:Eigenfunction}
\end{eqnarray}
Identifying eigenfunctions from data and using these linear models for control is a high-priority future direction.
\begin{figure}[b!]
\begin{center}
\vspace{-.25in}
\includegraphics[width=.95\textwidth]{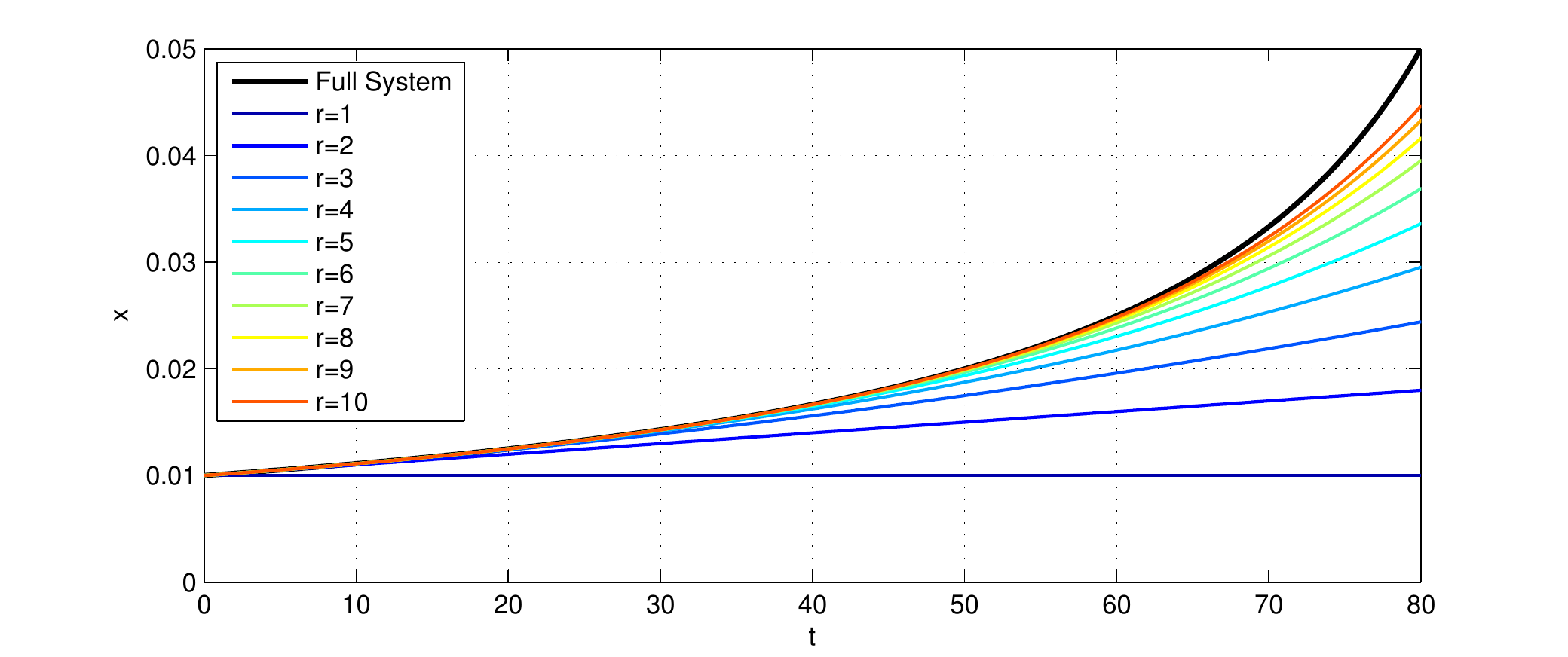}
\end{center}
\vspace{-.325in}
\caption{Illustration of Koopman linear system from Eq.~\eqref{Eq:CarlemanQuad} converging towards true solution as the rank of the truncation $r$ is increased.}\label{Fig:centermanifold_converge}
\end{figure}

\begin{lstlisting}[caption={Koopman linear system corresponding to Fig.~\ref{Fig:quad_manifold}.},label={Code:quad_manifold}]
clear all, close all, clc
%% System
mu = -.05;
lambda = -1;
A = [mu 0 0; 0 lambda -lambda; 0 0 2*mu];  % Koopman linear dynamics
[T,D] = eig(A);
slope_stab_man = T(3,3)/T(2,3);  % slope of stable subspace (green)

%% Integrate Koopman trajectories
y0A = [1.5; -1; 2.25];
y0B = [1; -1; 1];
y0C = [2; -1; 4];
tspan = 0:.01:1000;
[t,yA] = ode45(@(t,y)A*y,tspan,y0A);
[t,yB] = ode45(@(t,y)A*y,tspan,y0B);
[t,yC] = ode45(@(t,y)A*y,tspan,y0C);

%% Plot invariant surfaces
% Attracting manifold $y_2=y_1^2$ (red manifold)
[X,Z] = meshgrid(-2:.01:2,-1:.01:4);
Y = X.^2;
surf(X,Y,Z,'EdgeColor','None','FaceColor','r','FaceAlpha',.1)
hold on, grid on, view(-15,8), lighting gouraud

% Invariant set $y_3=y_1^2$ (blue manifold)
[X1,Y1] = meshgrid(-2:.01:2,-1:.01:4);
Z1 = X1.^2;
surf(X1,Y1,Z1,'EdgeColor','None','FaceColor','b','FaceAlpha',.1)

% Stable invariant subspace of Koopman linear system (green plane)
[X2,Y2]=meshgrid(-2:0.01:2,0:.01:4);
Z2 = slope_stab_man*Y2;  % for mu=-.2
surf(X2,Y2,Z2,'EdgeColor','None','FaceColor',[.3 .7 .3],'FaceAlpha',.7)

x = -2:.01:2;
% intersection of green and blue surfaces (below)
plot3(x,(1/slope_stab_man)*x.^2,x.^2,'-g','LineWidth',2) 
% intersection of red and blue surfaces (below)
plot3(x,x.^2,x.^2,'--r','LineWidth',2)  
plot3(x,x.^2,-1+0*x,'r--','LineWidth',2); 

%% Plot Koopman Trajectories (from lines 15-17)
plot3(yA(:,1),yA(:,2),-1+0*yA,'k-','LineWidth',1);
plot3(yB(:,1),yB(:,2),-1+0*yB,'k-','LineWidth',1);
plot3(yC(:,1),yC(:,2),-1+0*yC,'k-','LineWidth',1);
plot3(yA(:,1),yA(:,2),yA(:,3),'k','LineWidth',1.5)
plot3(yB(:,1),yB(:,2),yB(:,3),'k','LineWidth',1.5)
plot3(yC(:,1),yC(:,2),yC(:,3),'k','LineWidth',1.5)
plot3([0 0],[0 0],[0 -1],'ko','LineWidth',4)
set(gca,'ztick',[0 1 2 3 4 5])
axis([-4 4 -1 4 -1 4])
xlabel('y_1'), ylabel('y_2'), zlabel('y_3');
\end{lstlisting}


\begin{lstlisting}[caption={Koopman operator optimal control (KOOC) example corresponding to Fig.~\ref{Fig:lqr_nonlinear}.},label={Code:KOOC}]
clear all, close all, clc

mu = -.1;
lambda = 1;
tspan = 0:.01:50;
x0 = [-5; 5];

% LQR on linearized system
A = [-.1 0; 0 1];
B = [0; 1];
Q = eye(2);
R = 1;
C = lqr(A,B,Q,R);
vf = @(t,x) A*x + [0; -lambda*x(1)^2] - B*C*x;
[t,xLQR] = ode45(vf,tspan,x0);

% Koopman operator optimal control (KOOC); i.e., LQR on Koopman operator
A2 = [mu 0 0; 0 lambda -lambda; 0 0 2*mu];
B2 = [0; 1; 0];
Q2 = [1 0 0; 0 1 0; 0 0 0];
R = 1;
C2 = lqr(A2,B2,Q2,R);
% note that controller is nonlinear in the state 'x'
vf2 = @(t,x) A*x + [0; -lambda*x(1)^2] - B*C2(1:2)*x + [0; -C2(3)*x(1)^2];
[t,xKOOC] = ode45(vf2,tspan,x0);

%% Plot
figure(1)
subplot(1,3,1)
plot(xLQR(:,1),xLQR(:,2),'k','LineWidth',1.2);
hold on, grid on
plot(xKOOC(1:50:end,1),xKOOC(1:50:end,2),'r--','LineWidth',1.2);
xlabel('x_1'), ylabel('x_2')

subplot(1,3,2)
plot(tspan,xLQR,'k','LineWidth',1.2);
hold on, grid on
plot(tspan,xKOOC,'r--','LineWidth',1.2);
xlabel('t'), ylabel('x_k')
xlim([0 50])

JLQR = cumsum(xLQR(:,1).^2 + xLQR(:,2).^2 + (C*xLQR')'.^2)';
JKOOC = cumsum(xKOOC(:,1).^2 + xKOOC(:,2).^2 + (C*xKOOC')'.^2)';
subplot(1,3,3)
plot(tspan,JLQR,'k','LineWidth',1.2);
hold on, grid on
plot(tspan,JKOOC,'r--','LineWidth',1.2);
xlabel('t'), ylabel('J')
axis([0 50 0 500000])
legend('LQR','Koopman optimal control')
\end{lstlisting}

\end{document}